\documentclass[10pt]{amsart}

\usepackage{graphicx}
\usepackage{cite}
\usepackage[vcentermath]{youngtab}
\usepackage[margin=1.5in]{geometry}
\usepackage{mathtools}
\usepackage{ifthen}
\usepackage[unicode]{hyperref}

\theoremstyle{plain}
\newtheorem{theorem}{Theorem}[subsection]

\newtheorem{proposition}[theorem]{Proposition}

\theoremstyle{definition}
\newtheorem{definition}[theorem]{Definition}

\theoremstyle{remark}
\newtheorem{remark}[theorem]{Remark}
\newtheorem*{ack}{Acknowledgments}

\numberwithin{equation}{section}

\renewcommand{\C} {{\mathbb C}}                            
\newcommand{\R} {{\mathbb R}}                              
\newcommand{\Z} {{\mathbb Z}}                              
\newcommand{\N} {\Z_{\ge 0}}                               
\newcommand{\Pj} {{\mathbb P}}                             
\newcommand{\T} {{\mathbb T}}                              
\newcommand{\SL} {\mathsf{SL}} 

\newcommand{\suchthat} {\;\mid\;}

\newcommand*{\maps}{\colon}
\newcommand{\deftobe}{\mathrel{\coloneqq}}

\DeclareMathOperator{\Proj}{Proj}
\DeclareMathOperator{\Spec}{Spec}

\newcommand*{\gtp}[1]{\mathbf #1}

\newcommand*{\floor}[1]{\left\lfloor #1 \right\rfloor}

\renewcommand*{\epsilon}{\varepsilon}

\newcommand{\len}[1]{\operatorname{len}(#1)}

\newcommand*{\unilambda}{\9003\273}
\newcommand*{\unimu}{\9003\274}

\newcommand*{\ifpdfsyncstop}{
   \ifthenelse{
      \isundefined{\pdfsyncstop}
   }{
   }{
      \pdfsyncstop
   }
}
\newcommand*{\ifpdfsyncstart}{
   \ifthenelse{
      \isundefined{\pdfsyncstart}
   }{
   }{
      \pdfsyncstart
   }
}

\newlength{\Tij}
\newlength{\xij}
\title{%
   Degree bounds for type-A weight rings and Gelfand--Tsetlin
   semigroups
}%

\author{Benjamin J. Howard}
\address{Benjamin J. Howard \\ 
	 Mathematics Department \\
	 University of Michigan \\
	 Ann Arbor, MI 48109, USA}
\email{howardbj@umich.edu}

\author{%
   Tyrrell B. McAllister%
}
\address{%
   Tyrrell B. McAllister \\
   Department of Mathematics and Computer Science \\
   Eindhoven University of Technology \\
   PO Box 513 \\
   5600 MB Eindhoven \\
   The Netherlands%
}
\email{%
   tmcallis@win.tue.nl%
}
\thanks{%
   Second author supported by the Netherlands Organisation for
   Scientific Research (NWO) Mathematics Cluster DIAMANT%
}



\date{Draft}

\begin{document}

\begin{abstract}
   A weight ring in type A is the coordinate ring of the GIT
   quotient of the variety of flags in $\C^n$ modulo a twisted
   action of the maximal torus in $\SL(n,\C)$.  We show that any
   weight ring in type A is generated by elements of degree
   strictly less than the Krull dimension, which is at worst
   $O(n^2)$.  On the other hand, we show that the associated
   semigroup of Gelfand--Tsetlin patterns can have an essential
   generator of degree exponential in $n$.
\end{abstract}

\maketitle

\tableofcontents

\baselineskip=15pt

\section{Introduction}
   
Given a pair $\lambda, \mu$ of weights for $\SL_{n}(\C)$ with
$\lambda$ dominant, let $V_\lambda[\mu]$ denote the
$\mu$-isotropic component of the irreducible representation
$V_\lambda$ with highest weight $\lambda$.  The \emph{weight ring}
$R(\lambda,\mu)$ is the graded ring $\bigoplus_{N=0}^\infty V_{N
\lambda}[N \mu]$; it is the projective coordinate ring of the GIT
quotient of the flag variety modulo the $\mu$--twisted action of
the maximal torus $\T$ in $\SL_{n}(\C)$.  We define the
\emph{weight variety} $W(\lambda,\mu)$ as
\begin{equation*}
   W(\lambda,\mu) \deftobe \Proj R(\lambda,\mu).
\end{equation*}

\begin{remark}
   Weight varieties for arbitrary reductive Lie groups (not just
   those of type~A) were studied by A.~Knutson in his Ph.D. thesis
   \cite{Knutson}.  Knutson also studied the symplectic geometry
   of these spaces.
\end{remark}

Our first theorem (Theorem \ref{thm:KrullDimensionBound} below) is
that $R(\lambda,\mu)$ is generated in degree strictly less than
the Krull dimension of $R(\lambda,\mu)$, provided that the
degree-one piece $V_\lambda[\mu]$ is nonzero.  The basic idea
behind the proof is to show that the degree-one piece contains a
system of parameters, and that the $a$-invariant of
$R(\lambda,\mu)$ is negative.  (The $a$-invariant is the degree of
the Hilbert series, which is a rational function.)  The theorem
then follows from the fact that $R(\lambda,\mu)$ is
Cohen--Macaulay.
\begin{remark}
   The $a$-invariant is negative in all types; however, for types
   other than type~A, the degree-one piece sometimes fails to
   contain a system of parameters.  This condition is equivalent
   to the condition that all semistable flags lie in the supports
   of degree-one $\T$-invariants; see \cite{Howard} for a
   counterexample in $G = \mathrm{SO}_5(\C)$.
\end{remark}     
                       
It is well known (cf.  \cite{GonciuleaLakshmibai, KoganMiller,
FothHu}) that $R(\lambda,\mu)$ has a flat degeneration to the
semigroup algebra $R'(\lambda,\mu)$ of the semigroup of
Gelfand--Tsetlin patterns associated to semistandard tableaux of
shape $m\lambda$ and content $m\mu$ for $m \geq 0$.  In
particular, the ring $R'(\lambda,\mu)$ is the graded ring
associated to a filtration of $R(\lambda,\mu)$ by natural numbers.
Generators for $R'(\lambda,\mu)$ can be lifted to generators of
$R(\lambda,\mu)$, so one might hope that $R'(\lambda,\mu)$ is
relatively simple.  Unfortunately, we find pairs $\lambda,\mu$ for
which $R'(\lambda,\mu)$ has essential generators of degree
exponential in $n$.
      
Our second main result (Theorem \ref{thm:GTDegreeBound} below) is
that, in the case where $n = 3k$ is a multiple of $3$, the
semigroup algebra $R'(k\varpi_3,0)$ has an essential generator of
degree approximately $(\sqrt{2})^n$.  This is in striking contrast
to the lower bound of $2n-8$ for $R(k\varpi_3,0)$ that follows
from our first theorem (since the Krull dimension of $R(k
\varpi_3,0)$ is $2n-7$).  This case is particularly interesting,
because, via the Gelfand--MacPherson correspondence, it is the
moduli space of $n$-tuples of points in the projective plane.
This is a remarkable example of how a semigroup algebra produced
by a promising toric degeneration can fail to serve as an
effective proxy for the original ring.
   
Our motivation for studying the semigroup of Gelfand--Tsetlin
patterns was to imitate the method of \cite{HMSV}, which studied
the case of $n$ points on the projective line.  Here one takes
$\lambda$ to be a multiple of the second fundamental weight
$\varpi_{2}$.  It was shown in \cite{HMSV} that the associated
semigroup of Gelfand--Tsetlin patterns is generated in degree
$\leq 2$.  We had hoped to use the same method in the case of $n$
points in the projective plane, but the second theorem indicates
why this is not the right approach.  However, there might still be
another toric degeneration that yields a bound better than the one
in Theorem \ref{thm:KrullDimensionBound}, perhaps among those
discovered by Caldero \cite{Caldero}.
        
\begin{ack}
   We thank Harm Derksen, Ionut Ciocan-Fontanine, and Mircea
   Mustata for their invaluable advice.
\end{ack}

\section{%
   \texorpdfstring{%
      A description of the weight ring $R(\lambda, \mu)$%
   }{%
      Constructing the weight ring R(\unilambda, \unimu)%
   }%
}
\label{sec:WeightRingAsGenByProdOfMinors}

In this section, we give an explicit description of $R(\lambda,
\mu)$.  Let $n\ge 2$, and let $B$ denote the Borel subgroup of
$\SL_{n}(\C)$ consisting of the upper-triangular matrices in
$\SL_{n}(\C)$.  Fix a nontrivial dominant weight $\lambda$ of
$\SL_{n}(\C)$.  We represent $\lambda$ as a \emph{partition},
\emph{i.e.}, as a weakly decreasing sequence $(\lambda_{1},
\dotsc, \lambda_{n})$ of nonnegative integers, with $\lambda_{1}
\ge 1$ and $\lambda_n = 0$.  Let $\mu$ be a weight of
$\SL_{n}(\C)$ such that $V_{\lambda}[\mu]$ is nonzero.  Thus $\mu$
may be expressed as a sequence $(\mu_{1}, \dotsc, \mu_{n})$ of
nonnegative integers such that $\sum_{i=1}^{n} \mu_{i}
=\sum_{i=1}^{n} \lambda_{i}$.

By the Borel--Weil construction, the irreducible representation
with highest weight $\lambda$ is the finite-dimensional vector
space
\begin{equation*}
   V_\lambda
   =  \big\{
         \text{%
            holomorphic $f \maps \SL_{n}(\C) \to \C$%
         }%
         \suchthat
         \text{%
            $f(gb) = e^{\lambda}(b)f(g)$ for all $g \in \SL_{n}(\C)$
            and $b \in B$%
         }%
      \big\},
\end{equation*}
where $e^{\lambda}(b) \deftobe \prod_{i=1}^n b_{ii}^{\lambda_i}$
for $b = (b_{ij})_{1 \leq i,j \leq n} \in B$.  The action of
$\SL_{n}(\C)$ on $V_\lambda$ is given by $(g \cdot f)(h) =
f(g^{-1}h)$ for $g, h \in \SL_{n}(\C)$.  We define
$$R(\lambda) \deftobe \bigoplus_{N=0}^\infty V_{N \lambda}.$$
Multiplication in $R(\lambda)$ is the usual multiplication of
functions $\SL_n(\C) \to \C$.

\begin{remark}
   As we will review in the next section, the dominant weight
   $\lambda$ determines a line bundle $L_\lambda \to \SL_n(\C)/B$
   such that the space $\Gamma(\SL_n(\C)/B,L_\lambda)$ of sections
   is isomorphic to $V_{\lambda}$.  The multiplication in
   $R(\lambda)$ coincides with multiplication of the corresponding
   sections.  The ring $R(\lambda)$ is the coordinate ring of the
   partial flag variety.
\end{remark}

We define
\begin{equation*}
   V_\lambda[\mu] \deftobe \{f \in V_\lambda \mid \text{$f(tg) =
   e^\mu(t)f(g)$ for all $t \in \T$, $g \in \SL_n(\C)$}\}.
\end{equation*}
We now define a $\mu$-twisted action of $\T \subset \SL_n(\C)$ on
$R(\lambda)$.  For $f \in V_{N \lambda} = R(\lambda)_N$ of degree
$N$, the action of $t \in \T$ on $f$ is given by
\begin{equation*}
   (t \cdot f)(g) \deftobe e^{N\mu}(t) f(t^{-1}g).
\end{equation*}
Relative to this twisted action, the $\T$-invariant subring of
$R(\lambda)$ is exactly
\begin{equation*}
   R(\lambda,\mu) \deftobe \bigoplus_{N=0}^\infty
   V_{N\lambda}[N\mu].
\end{equation*}
 
\begin{remark}
   There is a unique $\SL_n(\C)$-linearization of $L_\lambda$.
   This defines a canonical $\T$-linearization of $L_\lambda$ by
   restriction $T \hookrightarrow \SL_n(\C)$.  The above action of
   $\T$ coincides with the canonical $\T$-linearization twisted by
   $\mu$.
\end{remark}

A fundamental fact from the representation theory of $\SL_{n}(\C)$
is that $V_{\lambda}$ has a basis indexed by semistandard tableaux
of shape $\lambda$.  A \emph{Young diagram} of shape $\lambda$ is
a left-justified arrangement of $\lambda_1 + \dotsb + \lambda_n$
boxes with $\lambda_{i}$ boxes in the $i$th row.  For example, if
$\lambda = (3,3,2,1,1,0)$, then the Young diagram of shape
$\lambda$ is
\begin{equation*}
   \yng(3,3,2,1,1).
\end{equation*}
A \emph{semistandard tableaux} of shape $\lambda$ is a filling of
each box in a Young diagram of shape $\lambda$ with a number from
$1$ through $n$ such that the rows are weakly increasing and the
columns are strictly increasing.  For example, if $\lambda =
(3,3,2,1,1,0)$ then
\begin{equation*}
   \tau = \young(115,246,35,5,6)
\end{equation*}
is a semistandard tableau of shape $\lambda$.

Such a tableau $\tau$ determines a basis vector $b_\tau \in
V_\lambda$ as follows.  Write $\len{I}$ for the length of a column
$I$ of $\tau$.  We identify $I$ with the $\len{I}$-tuple of its
entries, read from top to bottom.  If $I = (i_{1}, \dotsc,
i_{\len{I}})$, let $\det_{I} \maps \SL_{n}(\C) \to \C$ be the
function that returns the determinant of the $\len{I} \times
\len{I}$ submatrix consisting of rows $i_{1}, \dotsc, i_{\len{I}}$
and columns $1, 2, \dotsc, \len{I}$.  The basis vector $b_{\tau}$
is then defined by
\begin{equation*}
   b_{\tau} \deftobe \prod_{\text{columns $I$ of $\tau$}}
   \det\nolimits_{I}.
\end{equation*}
Hence, in the example above, $b_\tau = \det_{1,2,3,5,6}
\det_{1,4,5} \det_{5,6}$.

We can also describe the $\T$-isotropic subspace
$V_{\lambda}[\mu]$ in terms of semistandard tableaux.  The
\emph{content} of a tableau $\tau$ is $\mu = (\mu_1,\ldots,\mu_n)$
if $\mu_i$ is the number of boxes in $\tau$ that contain the
number $i$ for $1 \le i \le n$.  The subspace $V_\lambda[\mu]$ is
the span of the $b_\tau$ such that $\tau$ has shape $\lambda$ and
content $\mu$.

\section{%
   \texorpdfstring{%
      The first theorem: generators of $R(\lambda,\mu)$%
   }{%
      The first theorem: generators of R(\unilambda, \unimu)%
   }%
}

We will derive an upper bound on the degree in which
$R(\lambda,\mu)$ is generated.  We generally follow the method of
\cite{Popov} (also explained in \cite{DerksenKraft}).  The idea is
to find a homogeneous system of parameters, together with an upper
bound on the $a$-invariant of the ring; this yields an upper bound
on the degree of a generating set.

We begin by referencing a result of \cite{Howard} and showing why
this implies the existence of a system of parameters in degree
one.  First we must introduce the notion of semistability.  Recall
that $B$ is the Borel subgroup of $\SL_n(\C)$, and $\SL_n(\C)/B$
is the flag variety.  A function $f \in V_\lambda$ defines a
section of the line bundle $L_\lambda \deftobe \SL_n(\C) \times_B
\C \to \SL_n(\C)/B$, where $\SL_n(\C) \times_B \C$ denotes the
quotient of $\SL_n(\C) \times \C$ by the equivalence relation
$(gb,e^\lambda(b)z) \sim (g,z)$.  The projection $L_\lambda \to
\SL_n(\C)/B$ is given by sending the equivalence class of $(g,z)$
to $gB$.  We define the $\mu$-twisted linearization of $\T$ on
$L_\lambda$ by $t \cdot (g,z) \deftobe (t^{-1}g, e^\mu(t)z)$.
Given $f \in V_\lambda$, we define a global section $s_f$ of
$L_\lambda$ by $s_f(gB) = (g,f(g))$.  The map $f \mapsto s_f$ is
an isomorphism $V_\lambda \cong \Gamma(\SL_n(\C)/B, L_\lambda)$.
The $\mu$-twisted torus action on $L_\lambda$ defines an action on
global sections, which coincides with the $\mu$-twisted action on
$V_\lambda$ that we earlier defined.
  
A flag $gB$ is \emph{semistable} if there is a positive integer
$N$ and a $\T$-invariant global section $s$ of $L_{N \lambda}$
such that $s(gB) \neq 0$.  That is, $gB$ is semistable if and only
if there is an $N > 0$ and an $f \in V_{N \lambda}[N \mu]$ such
that $f(g) \neq 0$.  It was shown in \cite{Howard} that we may
take $N = 1$.  That is, $gB$ is semistable if and only if there
exists an $f \in V_\lambda[\mu]$ such that $f(g) \neq 0$.  We
shall use this fact to show that there is a system of parameters
within $V_\lambda[\mu]$ for $R(\lambda,\mu)$.
  
We now recall some basic facts from commutative algebra.  Our main
references on Cohen--Macaulay rings and modules are
\cite{BrunsHerzog, MillerSturmfels}.  Let $k$ be an algebraically
closed field.  Suppose that $A$ is a $\N$-graded
finitely-generated $k$-algebra with $A_0 = k$.  Let $\mathfrak{m}$
denote the graded ideal generated by the positive degree
homogeneous elements of $A$.  Then $\mathfrak{m}$ is the unique
graded ideal such that all other graded ideals are contained
within it.  A \emph{homogeneous system of parameters} for $A$ is a
set of homogeneous elements $x_1,\ldots,x_s$ such that $s$ is the
Krull-dimension of $A$ and the ideal $(x_1,\ldots,x_s)$ is
$\mathfrak{m}$-primary.  By \cite[Theorem 1.5.17]{BrunsHerzog} we
have that $x_1,\ldots,x_s$ is a homogeneous system of parameters
if and only if $A$ is an integral extension of the subalgebra
$k[x_1,\ldots,x_s]$, and that this is the case if and only if $A$
is a finitely-generated $k[x_1,\ldots,x_s]$-module.

Let the \emph{null cone} $\mathcal{N}$ be the subvariety of points
in $\Spec R(\lambda)$ at which all positive-degree homogeneous
elements of $R(\lambda,\mu)$ vanish.  The result of \cite{Howard}
translates into the following:
\begin{proposition}
   The elements of $V_\lambda[\mu]$ suffice to cut out the null
   cone set-theoretically.  That is, $\mathcal{N}$ is exactly the
   set of points at which all elements of $V_\lambda[\mu]$ vanish.
\end{proposition}
Now suppose that $I \subset R(\lambda)$ is the ideal of elements
vanishing on the null cone.  By the above proposition, $I$ is the
radical closure in $R(\lambda)$ of the ideal generated by
$V_\lambda[\mu] \subset R(\lambda)$.  Recall that $R(\lambda,
\mu)$ is the ring of polynomials in $R(\lambda)$ that are
invariant under the $\mu$-twisted action of $\T$.  Since $\T$ is
linearly reductive, there is a canonical projection $\pi:
R(\lambda) \to R(\lambda,\mu)$, called the \emph{Reynolds
operator}, which is $R(\lambda,\mu)$-linear.  Following Hilbert
(cf.  \cite[Prop.\ 3.1]{DerksenKraft}), we have the following
result.
\begin{proposition}
   The invariant ring $R(\lambda,\mu)$ is a finitely-generated
   module over the subalgebra generated by $V_\lambda[\mu] \subset
   R(\lambda,\mu)$.
\end{proposition}

\begin{proof}
   Let $J$ and $S$ be the ideal and subalgebra, respectively,
   generated by $V_\lambda[\mu]$ in $R(\lambda)$.  Then, since $I
   = \mathrm{Rad}(J)$, we have that $I^m \subset J$ for some $m >
   0$.  Since $\T$ is linearly reductive, the invariant ring
   $R(\lambda,\mu)$ is finitely generated.  Thus, there exist
   homogeneous $y_1,\ldots,y_t \in R(\lambda,\mu)$ such that
   $y_1,\ldots,y_t$ generate $R(\lambda,\mu)$.  Suppose that
   $h_1,\ldots,h_\ell$ span $V_\lambda[\mu]$.  We have that each
   $y_i^m$ belongs to the ideal $J$, and so $y_i^m =
   \sum_{j=1}^\ell f_j h_j$ for some homogeneous $f_j \in
   R(\lambda)$.  Now we apply the Reynolds's operator $\pi$ to
   obtain $y_i^m = \sum_{j=1}^\ell \pi(f_j) h_j$.  Each
   coefficient $\pi(f_j)$ is a homogeneous invariant of degree
   less than $y_i^m$.  It follows that $R(\lambda,\mu)$ is
   generated as an $S$-module by monomials $\mathbf{m} =
   \prod_{i=1}^t y_i^{e_i}$, where each $e_i < m$.  There are only
   a finite number of such monomials, proving the claim.
\end{proof}

We now have the following: (cf. Proposition 3.2 of \cite{DerksenKraft})

\begin{proposition}
   The degree-one piece $V_\lambda[\mu]$ of $R(\lambda,\mu)$
   contains a homogeneous system of parameters.
\end{proposition}

\begin{proof}
   We already know from the previous proposition that
   $R(\lambda,\mu)$ is a finitely-generated module over the
   subalgebra generated by its degree-one piece $V_\lambda[\mu]$.
   Let $s$ be the Krull dimension of $R(\lambda, \mu)$.  It is
   easy to show that $s$ generic elements $f_1, \dotsc, f_s$ in
   $V_\lambda[\mu]$ are algebraically independent and that
   $R(\lambda,\mu)$ is an integral extension of $\C[f_1, \dotsc,
   f_s]$.
\end{proof}

\begin{proposition}  
   The weight ring $R(\lambda,\mu)$ is Cohen--Macaulay.
\end{proposition}

\begin{proof}
   The argument is the same as that given in \cite[Corollary
   14.25]{MillerSturmfels} to show that $R(\lambda)$ is
   Cohen--Macaulay.  We know that $R(\lambda,\mu)$ has a Gr\"obner
   degeneration to a semigroup algebra of Gelfand--Tsetlin
   patterns (see Section \ref{sec:GensOfToricRing}).  Such
   semigroup algebras are the invariant subrings of polynomial
   rings by the action of a torus (cf.  \cite{Dolgachev}).  By the
   theorem of Hochster \cite{Hochster}, the subring of torus
   invariants in a polynomial ring is Cohen--Macaulay.  A general
   principle regarding Gr\"obner degenerations is that any good
   property of the special fiber is shared by the general fiber.
   This is true in particular for the Cohen--Macaulay property
   \cite[Corollary 8.31]{MillerSturmfels}.
\end{proof}     

\begin{proposition}
   If $f_1,\ldots,f_s$ is a homogeneous system of parameters for
   $R(\lambda,\mu)$, then $R(\lambda,\mu)$ is a free
   $\C[f_1,\ldots,f_s]$-module.
\end{proposition}

\begin{proof}
   Since $f_1,\ldots,f_s$ are algebraically independent, the ring
   $\C[f_1,\ldots,f_s]$ is regular, and so this proposition
   follows from \cite[Proposition 2.2.11]{BrunsHerzog}.
\end{proof}

For a graded module $M$, let $H(M;t) \deftobe \sum_{d=0}^\infty
\dim(M_d) \, t^d$ denote the \emph{Hilbert series} of $M$.  It is
well known that, if $M$ is finitely generated, then $H(M;t)$ is a
rational function in $t$.  Let $a(M)$ be the degree of $H(M;t)$ as
a rational function.  The number $a(M)$ is called the
\emph{$a$-invariant} of $M$.

Fix a homogeneous system of parameters $f_1,\ldots,f_s \in
V_\lambda[\mu]$ for $R(\lambda,\mu)$.  Let $S =
\C[f_1,\ldots,f_s]$ be the subalgebra generated by the $f_i$.  For
brevity of notation, we will write $R \deftobe R(\lambda,\mu)$.
Let $\mathbf{f}$ denote the $s$-tuple $f_1,\ldots,f_s$.  By
Theorems 13.37(5) and 13.37(6) of \cite{MillerSturmfels}, $R$ is a
free $S$-module, and
\begin{equation*}
   H(R/\mathbf{f}R;t) = H(R;t) (1-t)^s. 
\end{equation*}
But we can easily compute $H(R/\mathbf{f}R;t)$.  Suppose that $R =
Sy_1 \oplus \cdots \oplus Sy_m$.  Let $k \deftobe \max_j (\deg
y_j)$.  Now, $H(R/\mathbf{f}R;t)$ is the polynomial $p(t) =
\sum_{i=0}^k h_d t^d$, where $h_d$ is the number of $y_j$ such
that $\deg y_j = d$.  Therefore, we have proved the following.
   
\begin{proposition}\label{weakBound}
   The ring $R(\lambda,\mu)$ is generated in degree $\leq k = \dim
   R(\lambda,\mu) + a(R(\lambda,\mu))$.
\end{proposition}

\begin{proposition}\label{strongBound}
   The $a$-invariant of $R(\lambda,\mu)$ is negative.
\end{proposition}

\begin{proof}
   Let $R \deftobe R(\lambda,\mu)$.  The dimension of the $d$-th
   graded piece $R_d$ of $R$ is equal to the number of
   semistandard tableaux of shape $d\lambda$ and content $d\mu$;
   this coincides with the number of integer lattice points in the
   $d$-th dilate of the rational polytope $GT(\lambda,\mu)$ (see
   Definition \ref{def:GTpolytope} below).  As a result of the
   theory of lattice point enumeration for rational polytopes
   (see, \emph{e.g.}, \cite[Chapter 4]{Stanley}), we may conclude
   that the Hilbert series $H(R;t) = \sum_{d=0}^\infty f(d) t^d$
   is a rational function of negative degree.
\end{proof}
    
\begin{remark}
   In fact, in all types, given a pair of weights $\lambda,\mu$
   with $\lambda$ dominant, the dimension of the $d\mu$-weight
   space in the irreducible representation $V_{d\lambda}$ with
   highest weight $d\lambda$ equals the number of integer lattice
   points in the $d$-th dilate of a certain polytope (for example
   the string polytope associated with the reduced word for the
   longest Weyl element).  And so, in all types, the $a$-invariant
   of weight rings is negative.
\end{remark}
    
The above propositions imply our first theorem:

\begin{theorem}
\label{thm:KrullDimensionBound}
   The algebra $R(\lambda,\mu)$ is generated in degree strictly
   less than the Krull dimension of $R(\lambda,\mu)$.
\end{theorem}

\begin{proof}
   This follows immediately from Proposition \ref{weakBound} and
   Proposition \ref{strongBound}.
\end{proof}

Finally, we point out that the Krull dimension of $R(\lambda,\mu)$
is one more than the dimension of the GIT quotient of the flag
variety by $\T$.  This is at most the dimension of the flag
variety itself, which is $n(n-1)/2$.  In the case of $n$ points in
projective space $\Pj^{m-1}$, where $\lambda$ is a multiple of the
$m$-th fundamental weight $\varpi_m$ for $\SL_{n}(\C)$, the Krull
dimension of $R(\lambda,\mu)$ is at most $n(m-1)-(m^2-1)+1$.

\section{The toric degeneration to Gelfand--Tsetlin patterns}
\label{sec:GensOfToricRing}

A \emph{Gelfand--Tsetlin pattern}, or GT~pattern, is a triangular
array $\gtp{x} = (x_{ij})_{1 \le i \le j \le n}$ of real numbers
satisfying the \emph{interlacing inequalities} $x_{i,j+1} \ge
x_{ij} \ge x_{i+1,j+1}$.  We express $\gtp{x}$ as a triangular
array by arranging the entries as follows:
\begin{equation*}
   \begin{matrix}
      x_{1n} &        & x_{2n} &        & x_{3n} & \multicolumn{3}{c}{\cdots}
                                                                            & x_{nn} \\
             &        &        &        &        &        &        &        &        \\
             & \hdotsfor{7}                                                          \\
             &        &        &        &        &        &        &        &        \\
             &        & x_{13} &        & x_{23} &        & x_{33} &        &        \\
             &        &        &        &        &        &        &        &        \\
             & \phantom{x_{02}}
                      &        & x_{12} &        & x_{22} &        & \phantom{x_{32}}
                                                                            &        \\
             &        &        &        &        &        &        &        &        \\
             &        &        &        & x_{11} &        &        &        &
   \end{matrix}
\end{equation*}

Given a semistandard tableaux $\tau$ with entries from $1$ through
$n$, let $\tau(j)$ be the tableau obtained from $\tau$ by deleting
all boxes containing indices strictly larger than $j$.  Hence,
$\tau(n) = \tau$.  Let $\lambda(j)$ denote the shape of $\tau(j)$.
One obtains an integral GT~pattern $\gtp{x}(\tau) =
(x(\tau)_{ij})_{1 \le i \le j \le n}$ by letting $x(\tau)_{ij} =
\lambda(j)_i$.  If $\tau$ has shape $\lambda$ and content $\mu =
(\mu_1,\ldots,\mu_n)$, then the resulting GT~pattern
$\gtp{x}(\tau)$ has top row $\lambda$, and, for $1 \leq j \leq n$,
\begin{equation*}
   \sum_{i=1}^j x(\tau)_{ij} = \mu_1 + \cdots + \mu_j.
\end{equation*}
We denote this assignment by $\Phi \maps \tau \to \gtp{x}(\tau)$.
It is easy to see that it is a bijection from semistandard
tableaux of shape $\lambda$ and content $\mu$ to integral
GT~patterns with top row $\lambda$ and row sums equal to the
partial sums of $\mu$.  The GT~patterns with a fixed top row and
fixed row sums constitute a rational polytope.

\begin{definition}\label{def:GTpolytope}
   The \emph{GT~polytope} $GT(\lambda,\mu)$ is the set of
   \emph{real} GT patterns $(x_{ij})_{1 \le i \le j \le n}$ with
   top row $\lambda$ and with row sums $\sum_{i=1}^j x_{ij} =
   \mu_1 + \cdots + \mu_j$ for $1 \le j \le n$.
\end{definition}

Let $S(\lambda, \mu)$ denote the graded semigroup of integer
GT~patterns (under addition) that lie in $GT(N\lambda,N\mu)$ for
some nonnegative integer $N$.  Gonciulea and Lakshmibai have
described a Gr\"obner degeneration of the ring $R(\lambda) =
\bigoplus_{N=0}^\infty V_{N \lambda}$ to a semigroup algebra
$R'(\lambda)$ as the special fiber~\cite{GonciuleaLakshmibai}.  It
was shown in \cite{KoganMiller} (and also in \cite[Corollary
14.24]{MillerSturmfels}) that this semigroup is isomorphic to the
semigroup of integral GT~patterns with top row $N\lambda$ for some
nonnegative integer $N$.  This construction also applies to the
subring $R(\lambda,\mu)$ by restricting to $\T$-invariants, as we
now describe.  See~\cite{FothHu} for details.
 
The resulting degenerated ring $R'(\lambda,\mu)$ has the same
underlying graded vector space as $R(\lambda,\mu)$.  The
semistandard tableaux of shape $N\lambda$ and content $N\mu$, $N >
0$, index a basis for $R'(\lambda,\mu)_{N}$.  Let $b'_\tau \in
R'(\lambda,\mu)$ denote the basis element indexed by $\tau$.  The
basis element $b'_\tau$ is the leading term of $b_\tau \in
R(\lambda,\mu)$ for a certain filtration of $R(\lambda,\mu)$ (see
\cite{GonciuleaLakshmibai} and \cite{KoganMiller, FothHu}).  The
filtration has the special property that, if $\tau_1,\tau_2$ are
any two semistandard tableaux, and if $b_{\tau_1}b_{\tau_2} =
\sum_\tau c^\tau b_\tau$, where the sum is over semistandard
tableaux, then the term $b_{\Phi^{-1}(\gtp{x}(\tau_1) +
\gtp{x}(\tau_2))}$ appears on the right-hand side with coefficient
equal to $1$.  Furthermore, all other terms $c^\tau b_\tau$ have
strictly smaller filtration level.  Thus, in $R'(\lambda,\mu)$,
the multiplication rule becomes
\begin{equation*}
   b'_{\tau_1}b'_{\tau_2} =
   b'_{\Phi^{-1}(\gtp{x}(\tau_1)+\gtp{x}(\tau_2))}.
\end{equation*}
Therefore $R'(\lambda,\mu)$ is isomorphic to $\C[S(\lambda,
\mu)]$, the semigroup algebra of GT~patterns under addition of
patterns.

Given an $m$-tuple of rational numbers $q_1,\ldots,q_m$, define
$\text{den}(q_1,\ldots,q_m)$ to be the least positive integer $N$
such that $Nq_i \in \Z$ for each $i$, $1 \leq i \leq m$.  We call
this the \emph{denominator} of the $m$-tuple.  Now, if some vertex
$\gtp{x}$ of the polytope $GT(\lambda,\mu)$ has denominator $N >
1$, then the integer point $N \gtp{x}$ is an essential generator
of the semigroup $\C[S(\lambda, \mu)]$, since $\gtp{x}$ cannot be
written as a sum of other integral patterns in $S(\lambda, \mu)$.
In the next section we show the existence of such a vertex with
large denominator for the case where $n$ is a multiple of $3$,
$\lambda = \tfrac{n}{3} \varpi_3$, and $\mu = (1,1,\ldots,1)$.

\section{%
   \texorpdfstring{%
      The second theorem: $3k$ points on $\Pj^2$ and a nasty
      GT~pattern%
   }{%
      The second theorem: 3k points on P\texttwosuperior%
   }
}
\label{sec:ExpDegBound}

Suppose that $n = 3k$, where $k \ge 2$ is an integer.  Let
$\lambda = k\varpi_{3}$ be a multiple of the third fundamental
weight for $\SL_{n}(\C)$.  Thus, as a partition, $\lambda = (k, k,
k, 0, \dotsc, 0) \in \R^{3k}$.  Now let $\mu$ be the
``democratic'' weight dominated by $\lambda$.  That is, we
represent $\mu$ by the composition $(1, \dotsc, 1) \in \R^{3k}$.
With this choice of $\lambda$ and $\mu$, the projective variety
$\Proj R(\lambda, \mu)$ is the moduli space of equally weighted
$3k$-tuples of points in projective space $\Pj^2$ (see \cite{door}
for more details).

We now construct a GT~pattern that we claim will be a vertex of
$GT(\lambda, \mu)$.  Define the sequences $\{T^{(1)}_{j}\}$ and
$\{T^{(2)}_{j}\}$ by the coupled recurrence relations
\begin{equation}
\label{T1jT2jRecDef}
   \begin{aligned}
      T^{(1)}_{0} &= k & T^{(2)}_{0} &= k - 1/2 \\
      T^{(1)}_{j} &= T^{(2)}_{j-1} - 1 \quad (j \ge 1) \qquad \qquad
         & T^{(2)}_{j}
            &= \frac{1}{2}
               \left(
                   T^{(1)}_{j} + T^{(1)}_{j-1}
               \right) \quad (j \ge 1).
   \end{aligned}
\end{equation}
Solving this system of recurrence relations yields the closed-form
expressions
\begin{align}
   T^{(1)}_{j}
   &= k - \frac{2}{3}j + \frac{5}{9}
      \left(\frac{-1}{2}\right)^{j} - \frac{5}{9}
      \label{eq:T1jClosedForm} \\
   T^{(2)}_{j}
   &= k - \frac{2}{3}j - \frac{5}{18}
      \left(\frac{-1}{2}\right)^{j} - \frac{2}{9}.
      \label{eq:T2jClosedForm}
\end{align}

Let $N = k + \floor{k/2} - 2$.  We will construct a triangular
array $\gtp{x}$ by filling in the entries of $\gtp{x}$ in blocks
from the upper left to the lower right using the values
$T^{(1)}_{j}$ and $T^{(2)}_{j}$.  Begin by filling the entries in
the upper left of the triangular array as follows.
\begin{equation*}
   \begin{matrix}
      \makebox[\xij]{$x_{1n}$} &  & \makebox[\xij]{$x_{2n}$} &  & \makebox[\xij]{$x_{3n}$} &   \\
       &  &  &  &  &   \\
       & \makebox[\xij]{$x_{1,n-1}$} &  & \makebox[\xij]{$x_{2,n-1}$} &  & \makebox[\xij]{$x_{3,n-1}$}  \\
       &  &  &  &  &   \\
       &  & \makebox[\xij]{$x_{1,n-2}$} &  & \makebox[\xij]{$x_{2,n-2}$} &   \\
       &  &  &  &  &   \\
       &  &  & \makebox[\xij]{$x_{1,n-3}$} &  &
   \end{matrix}
   \qquad
   =
   \quad
   \begin{matrix}
      \makebox[\xij]{$k$} &  & \makebox[\xij]{$k$} &  & \makebox[\xij]{$k$} &   \\
       &  &  &  &  &   \\
       & \makebox[\xij]{$k$} &  & \makebox[\xij]{$k$} &  & \makebox[\xij]{$k-1$}  \\
       &  &  &  &  &   \\
       &  & \makebox[\xij]{$k$} &  & \makebox[\xij]{$k - \frac{1}{2}$} &   \\
       &  &  &  &  &   \\
       &  &  & \makebox[\xij]{$k$} &  &
   \end{matrix}
\end{equation*}
We then proceed from the upper left to the lower right of the
triangular array by filling in blocks of entries as follows.
For $1 \le j \le N-1$, let
\begin{equation*}
   \begin{matrix}
       &  & \makebox[\xij]{$x_{3, n - 2j}$} &   \\
       &  &  &   \\
       & \makebox[\xij]{$x_{2, n - 2j - 1}$} &  & \makebox[\xij]{$x_{3, n - 2j - 1}$}  \\
       &  &  &   \\
       \makebox[\xij]{$x_{1, n - 2j - 2}$} &  & \makebox[\xij]{$x_{2, n - 2j - 2}$} &   \\
       &  &  &   \\
       & \makebox[\xij]{$x_{1, n - 2j - 3}$} &  &
   \end{matrix}
   \qquad
   =
   \qquad
   \begin{matrix}
       &  & \makebox[\xij]{$T^{(1)}_{j}$} &   \\
       &  &  &   \\
       & \makebox[\xij]{$T^{(1)}_{j}$} &  & \makebox[\xij]{$T^{(1)}_{j}$}  \\
       &  &  &   \\
       \makebox[\xij]{$T^{(2)}_{j}$} &  & \makebox[\xij]{$T^{(1)}_{j}$} &   \\
       &  &  &   \\
       & \makebox[\xij]{$T^{(1)}_{j}$} &  &
   \end{matrix}
\end{equation*}
If $k$ is even, the final entries at the bottom of the array are
filled in as follows.
\begin{equation*}
   \begin{matrix}
       &  &              x_{34} &   \\
       &  &  &   \\
       &        x_{23} &   &      x_{33}  \\
       &  &  &   \\
       x_{12} &   &      x_{22} &   \\
       &  &  &   \\
       &        x_{11} &  &
   \end{matrix}
   \qquad
   =
   \qquad
   \begin{matrix}
       &  &  T^{(1)}_{N} &  \\
       &  &  &   \\
       & T^{(1)}_{N} &  & T^{(1)}_{N}  \\
       &  &  &   \\
       \makebox[\Tij]{$2 - T^{(1)}_{N}$} &  & T^{(1)}_{N} &  \\
       &  &  &   \\
       & 1 &  &
   \end{matrix}
\end{equation*}
On the other hand, if $k$ is odd, then the final entries
are filled in as follows:
\begin{equation*}
   \begin{matrix}
       &  & x_{35} &  &   \\
       &  &  &   \\
       & x_{24} &  & x_{34} &   \\
       &  &  &   \\
       x_{13} &  & x_{23} &  & x_{33}  \\
       &  &  &   \\
       & x_{12} &  & x_{22} &   \\
       &  &  &   \\
       &  & x_{11} &  &
   \end{matrix}
   \qquad
   =
   \qquad
   \begin{matrix}
       &  & T^{(1)}_{N} &  &  \\
       &  &  &   \\
       & T^{(1)}_{N} &  & T^{(1)}_{N} & \\
       &  &  &   \\
       T^{(1)}_{N} &  & T^{(1)}_{N} &  & \makebox[\Tij]{$3 - 2T^{(1)}_{N}$}  \\
       &  &  &   \\
       & T^{(1)}_{N} &  & \makebox[\Tij]{$2 - T^{(1)}_{N}$} & \\
       &  &  &   \\
       &  & 1 &  &
   \end{matrix}
\end{equation*}
All the remaining entries of the triangular array are assigned the
value 0.

\ifpdfsyncstop
\begin{figure}[tbp]
   \includegraphics{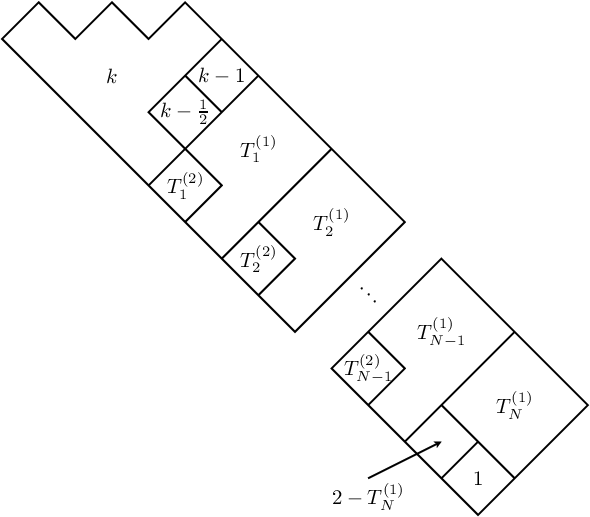}
   \caption{Tiling of $\gtp{x}$ when $k$ is even}
   \label{fig:Tiling}
\end{figure}
\ifpdfsyncstart

\begin{proposition}
\label{prop:ArrayIsVertex}
   The triangular array constructed above is a vertex of
   $GT(\lambda, \mu)$ with denominator $2^{N}$.
\end{proposition}

\begin{proof}
   To show that $\gtp{x} \in GT(\lambda, \mu)$, we first check
   that $\gtp{x}$ is a GT-pattern.  In this case, the interlacing
   inequalities to be verified are
   \begin{align*}
      \left.
         \begin{gathered}
            T^{(1)}_{j-1} > T^{(1)}_{j} > 0   \\
            T^{(1)}_{j-2} > T^{(2)}_{j-1} > T^{(1)}_{j-1}
         \end{gathered}
      \right\}
         & \quad \text{for $2 \le j \le N$,} \\
      \\
      \left.
         \begin{gathered}
            T^{(1)}_{N-1} > 2 - T^{(1)}_{N} > T^{(1)}_{N}  \\
            2 - T^{(1)}_{N} > 1 > T^{(1)}_{N}
         \end{gathered}
      \right\}
         & \quad \text{if $k$ is even,} \\
      \\
      \left.
         \begin{gathered}
            T^{(1)}_{N} > 3 - 2  T^{(1)}_{N} > 0  \\
            T^{(1)}_{N} > 2 - T^{(1)}_{N} > 3 -  2T^{(1)}_{N}  \\
            T^{(1)}_{N} > 1 >  2 - T^{(1)}_{N}
         \end{gathered}
      \right\}
         & \quad \text{if $k$ is odd.}
   \end{align*}

   These are all straightforward consequences of the closed-form
   expressions \eqref{eq:T1jClosedForm} and
   \eqref{eq:T2jClosedForm} for $T^{(1)}_{j}$ and $T^{(2)}_{j}$,
   respectively, so $\gtp{x}$ is a GT-pattern.  Thus, to show that
   $\gtp{x} \in GT(\lambda, \mu)$, we need only establish that the
   row-sums of $\gtp{x}$ are correct.  This amounts to showing
   that
   \begin{align}
      T^{(2)}_{j-1} + T^{(1)}_{j-1} + T^{(1)}_{j} & = 3k - 2j \notag \\
      T^{(1)}_{j-1} + 2 T^{(1)}_{j} & = 3k - 2j - 1
      \label{AltT1jDef}
   \end{align}
   for $2 \le j \le N$.  These equalities may be shown using
   induction and the recursive definition \eqref{T1jT2jRecDef} of
   $T^{(1)}_{j}$ and $T^{(2)}_{j}$.  It is clear from equation
   \eqref{AltT1jDef} that $T^{(1)}_{j}$ has denominator $2^{j}$
   when written as a reduced fraction.  Hence, $\gtp{x}$ has
   denominator $2^{N}$, as claimed.

   It remains only to show that $\gtp{x}$ is a vertex of
   $GT(\lambda, \mu)$.  We prove this by showing that, for any
   triangular array $\epsilon$, if $\gtp{x} \pm \epsilon \in
   GT(\lambda, \mu)$, then $\epsilon = 0$.  This is most easily
   seen by partitioning the entries of $\gtp{x}$ so that entries
   that are equal and adjacent are grouped together.  We call each
   group of entries in this partition a \emph{tile}.  See Figure
   \ref{fig:Tiling} for a depiction of the case when $k$ is even.
   Each tile is labeled with the value shared by the entries that
   it contains.

   Suppose that $\gtp{x} \pm \epsilon \in GT(\lambda, \mu)$.  Note
   that, after the addition of $\pm \epsilon$, the entries in each
   tile must still share a value, and the row-sums must be
   unchanged.  We prove inductively that the entries in each tile
   cannot have changed, proceeding from the upper left to the
   lower right.

   The entries in the tile labeled $T^{(1)}_{0} = k$ cannot have
   changed because the top row is fixed.  For the same reason, the
   $0$ entries in $\gtp{x}$ are also fixed.  Proceeding by
   induction, the entries in the tile labeled $T^{(1)}_{j}$ cannot
   have changed because there is a row on which this is the only
   tile besides the tile labeled $T^{(1)}_{j-1}$ and the tile of
   $0$s, which have already been fixed.  Hence, the entries in all
   the tiles labeled $T^{(1)}_{j}$, $0 \le j \le N$, are fixed
   under the addition of $\pm\epsilon$.  Finally, for $1 \le j \le
   N-1$, the tile labeled $T^{(2)}_{j}$ lies on a row in which the
   other entries, $T^{(1)}_{j}$, $T^{(1)}_{j+1}$, and 0, have been
   shown to be fixed, so the entry in this tile is also fixed
   under the addition $\pm \epsilon$.  Therefore, we conclude that
   $\epsilon = 0$, so that $\gtp{x}$ is a vertex, as claimed.
\end{proof}

The following theorem is an immediate consequence.
\begin{theorem}
\label{thm:GTDegreeBound}
   The Gelfand--Tsetlin algebra $R'(k \varpi_{3}, \mu)$ has
   essential generators of degree exceeding $2^{n/2-3}$.
\end{theorem}


\def\cprime{$'$}
\providecommand{\bysame}{\leavevmode\hbox to3em{\hrulefill}\thinspace}
\providecommand{\MR}{\relax\ifhmode\unskip\space\fi MR }
\providecommand{\MRhref}[2]{%
  \href{http://www.ams.org/mathscinet-getitem?mr=#1}{#2}
}
\providecommand{\href}[2]{#2}

\end{document}